\documentclass{amsart}      

\usepackage[psamsfonts]{amssymb}
\usepackage[psamsfonts]{eucal}
\usepackage{epsfig}
\usepackage{color}
\usepackage{graphicx}
\usepackage{url}

\title [An Investigation of the Collatz conjecture]{An Investigation of the Collatz conjecture
\newline
\\ \Small{John G. Koelzer$^\dagger$, Rockhurst University \\Daniel J. Welling}
\\ \today\ }
\thanks{$^\dagger$email: john.koelzer@rockhurst.edu}

\begin{document}             

\textcolor{blue}{
\maketitle }                 
\thispagestyle{empty}	     

\label{Abstract}
\noindent
\textcolor {blue}{\textbf{Abstract:} }  This paper explores special conditions on the starting value of a Collatz sequence which imply that the Collatz conjecture is true.  This is the result of the collaboration of a retired mathematics professor (Koelzer) and a retired physics professor (Welling).

\newpage
\label{Definition Collatz Conjecture}
\noindent
The Collatz conjecture was formulated by L. Collatz in 1937.  The conjecture concerns the definition of a sequence of positive integers by means of a simple algorithm.  The conjecture states that no matter what the starting value is, the sequence eventually equals 1.  For more information on the history of the Collatz problem see \cite{1} and \cite{2}.

\noindent
\textcolor {blue}{\textbf{The Collatz Conjecture:} } 
Let  $n$ be a positive integer.  If $n$ is even, divide it by 2 to get $n / 2$.  If $n$ is odd, multiply it by 3, add 1 and divide the result by 2 to obtain the number $(3n + 1)/2$.  Repeat the process indefinitely. The conjecture is that no matter what number you start with, you will always reach 1.   Our definition follows that of Terras in \cite{3}.  This is equivalent to the original Collatz algorithm but it results in a somewhat shorter sequence.  We will assume that our starting value is an odd integer; otherwise we can simply divide by 2 until the number is odd.

\noindent
\textcolor{blue}{\textbf{Example:} } The Collatz sequence for $n = 19$ is:  $19\to29\to44\to22\to11\to17\to26\to13\to20\to10\to5\to8\to4\to2\to1$. (14 steps)

In this paper we will present some results concerning the Collatz conjecture. We will assume that the starting integer $N $ in a Collatz sequence is an odd integer, since, if $N$ is even, it can be repeatedly be divided by $2$ until the resulting integer is odd.

The following theorem expresses the relation between an odd integer in a Collatz sequence and the next even integer in the sequence.  

\label{Theorem 1}
\noindent
\textcolor {blue}{\textbf{Theorem 1:} } Let $N$ be an odd integer in a Collatz sequence and define $E$ to be $N+1$.  The successive odd integers in the sequence starting with   $N$ are

\begin{center}
$3\left(\dfrac{E}{2}\right)-1,\  3^2\left(\dfrac{E}{2^2}\right)-1, \  3^3\left(\dfrac{E}{2^3}\right)-1,\  \cdots, \  3^{k-1}\left(\dfrac{E}{2^{k-1}}\right)-1,$
\end{center}

\noindent
where $k$  is the highest power of 2 contained as a factor in $E$.  The integer $3^k\left(\dfrac{E}{2^k}\right)-1$ is the first even integer after $N$. 

\noindent
\textbf {Proof:}  Let $m$ be an integer such that $1 \le m \le k-1.$ Then $\dfrac{E}{2^m}$ is an even number and hence $3^m\left(\dfrac{E}{2^m}\right)-1$ will be odd.  We will show that these integers are the successive odd entries in the Collatz sequence starting with $N$.

\noindent
\textbf{Case 1: } $m=1$.
The next term in the Collatz sequence after $N$ is  

\begin{center} 
$\dfrac{3N+1}{2}=\dfrac{3}{2}N+\dfrac{1}{2}=\dfrac{3}{2}\left(E-1\right)+\dfrac{1}{2}=\dfrac{3}{2}E-1$,
\end{center}

\noindent
 and so the theorem is proven for $m = 1$.
 
\noindent
\textbf{Case 2:} $2 \le m \le k-1.$
Since $3^{m-1}\left(\dfrac{E}{2^{m-1}}\right)-1$ is an odd integer, the next term in the Collatz sequence is $\dfrac{3\left(3^{m-1}\left(\dfrac{E}{2^{ m-1}}\right)-1\right)+1}{2}$, which simplifies algebraically to $3^m\left(\dfrac{E}{2^m}\right)-1$ and the theorem is proved in this case.Finally, since $k$ is the highest power of 2 contained as a factor in $E$, $\dfrac{E}{2^k}$ will be an odd integer.  This means that $3^k\left(\dfrac{E}{2^k}\right)-1$ is even and it is the first even integer in the Collatz sequence after $N$.   

Given an odd integer $N$, let $u$ be the number of steps to the next even integer, say $M$, and let $d$ be the number of steps from $M$ to the next odd integer $N'$.  By Theorem 1  $u$ is the largest power of 2 in $N+1$ and, since $M$ is even, $d$ is the largest power of 2 in $M$.  The relation between $N$ and $N'$ is given by the formula
\begin{center}
	$N' =\left [(\frac{3}{2})^{u}(N-1\right](\frac{1}{2})^{d}$
\end{center}

Here is a simple corollary to illustrate the theorem:

\label{Corollary 1}
\noindent
\textcolor{blue}{\textbf {Corollary 1:}} Suppose $k$ is a positive integer.  Let $N = 10^k - 1$ be an odd integer in a Collatz sequence.  Then the next even integer in the sequence will be $15^k - 1$.

\noindent
\textbf{Proof:} Using the notation from Theorem 1, $E = 10^k$ and the largest power of 2 in $E$ is $2^k$. So by the Theorem, the first even integer following $10^k-1$ is $3^k(\tfrac{E}{2^k})-1 =15^k-1$.

Another result following from Theorem 1 is shown below:

\label{Theorem 2}
\noindent
\textcolor {blue}{\textbf{Theorem 2:}} Let $k$, $l$ and $m$ be positive integers with $m$ odd and let $N = \left(\frac{2}{3}\right)^k\left(2^l m+1\right)-1$.  If $N$ is an odd integer then $N$ reaches $m$ in $l+k$ steps.

\noindent
\textbf{Proof:} Let $E$ be the even integer $N +1$.  Note that $E$ can be written as $\dfrac{2^k\left(2^l m+1\right)}{3^k}$.  We claim that $k$ is the highest power of 2 contained as a factor in $E$.  Suppose $p$ is the highest power; Since $2^k$ is a factor of $E$, we must have $k \le p$.

Furthermore, we can write $\dfrac{E}{2^p}$ as $\dfrac{2^k\left(\dfrac{2^l m+1}{3^k}\right)}{2^p}$.  Since $2^l m+1$ and $3^k$ are odd, their quotient is odd and so $2^p$ must divide $2^k$.  This implies $p \le k$ and we conclude $k = p$.

By Theorem 1, $N$ is followed by $k-1$ odd integers in the Collatz sequence and the even number $3^k\left(\dfrac{E}{2^k}\right)-1$.  This even number can be rewritten as $\left(\dfrac{3}{2}\right)^k\left[\dfrac{2^k(2^l m+1)}{3^k}\right]-1$, which algebraically reduces to $2^lm$.  This means that the next $l$ terms of the Collatz sequence are all products of $m$ multiplied by powers of 2 terminating at $m$ and the total number of steps from $N$ to $m$ is $l+k$.

\noindent
\textcolor{blue}{\textbf{Example:} }
$k=3, l=2, d=47\ \to N=(8/27)(4\cdot47+1)-1=55$. The Collatz sequence from 55 to 47 is: $55\to 83\to 125\to 188\to 94\to 47$. (5 steps)

\noindent
\textcolor{blue}{\textbf{Remarks: }} To apply Theorem 1, it is required to find integers $k$, $l$ and $m$ that make $N$ odd.  Also, it would be nice to show that there is an infinite number of values of $k$, $l$ and $m$ satisfying the hypotheses of the theorem, but that is not essential to the proof.

Here is a special case of Theorem 2:

\label{Corollary 2}
\noindent
\textcolor {blue}{\textbf{Corollary 2:} }Let $k$ and $l$ be positive integers and let $N = \left(\frac{2}{3}\right)^k\left(2^l +1\right)-1$.  If $N$ is an odd integer then $N$ reaches $1$ in $l+k$ steps.

\noindent
\textbf{Proof: } Let $m=1$ in Theorem 1.

\noindent
\textcolor{blue}{\textbf{Examples:}}\\
$k = 1, l = 5 \to N=\frac{2}{3}\left(2^5+1\right)-1 = 21.$
The Collatz sequence for 21 is:  $21\to 32\to16\to 8\to4\to2\to 1.$ (6 steps)

\noindent
$k = 2, l = 9 \to N= \left(\frac{2}{3}\right)^2\left(2^9+1\right)-1 = 227.$
The Collatz sequence for 227 is:  $227\to341\to512\to 256\to128\to 64\to 32\to16\to 8\to4\to2\to1$.  (11 steps)

\noindent
$k = 3, l = 9 \to N=\left(\frac{2}{3}\right)^3\left(2^9+1\right)-1 = 151.$ 
The Collatz sequence for 151 is:  $ 151\to227\to341\to512\to256\to128\to64\to \cdots \to4\to2\to1.$   (12 steps)

We prove a result that will be useful later in the paper.

\label{Lemma 1}
\noindent
\textcolor {blue}{\textbf{Lemma 1:} } Let $n$ be a positive integer. Then  $2^{(3^{n-1})} \equiv -1 \mod{3^n}$.

\noindent
\textbf{Proof: } We will use mathematical induction ~\cite{4}.

\noindent
\textbf{Case 1:} If $n=1$, $2^{(3^0)} = 2 \equiv -1 \mod{3^1}$

\noindent
\textbf{Case 2:}  Assume the lemma is true for $n=k$;  i.e., $2^{(3^{k-1})} \equiv -1 \mod{3^k}$.  We may write 
$2^{(3^{k-1})}  = 3^kq-1$ for some integer $q$. To show the equation is true for $n = k+1$, we  cube both sides of this equation using the binomial theorem:
\noindent
\begin{align}
2^{3^{k}}     &={(2^{3^{k-1}})}^3                           \notag \\
              &=\left(3^kq-1\right)^3                       \notag \\
              &=(3^kq)^3-3(3^kq)^2+3(3^kq)-1                  \notag \\
              &=q^33^{3k}-3q^23^{2k}+3(3^kq)-1                \notag \\
              &=q^33^{2k-1}3^{k+1}-q^23^k3^{k+1}+3^{k+1}q-1 \notag \\
              &=\left(q^33^{2k-1}-q^23^k+q\right)3^{k+1}-1  \notag 
\end{align}
\noindent
So $2^{3^{k}} -(-1)$ is a multiple of $3^{k+1}$ and the formula is verified for $n = k+1$.

\noindent
\textcolor{blue}{\textbf{Example:} }For $n=4$, $2^{(3^{n-1})}= 2^{27}=134,217,728 \equiv 80 \equiv -1 \mod{81}$.

\label{Theorem 3}
\noindent
\textcolor {blue}{\textbf{Theorem 3:}} Suppose $k$ is a positive integer and $l = 3^{k-1}\cdot r$, where $r$ is an odd integer.  If   $N = \left(\frac{2}{3}\right)^k\left(2^l+1\right)-1$  then the Collatz sequence beginning with $N$ terminates at 1 after $l+k$ steps.

\noindent
\textbf{Proof: }  By invoking Corollary 1, it suffices to show that $N$ is an odd integer under the given hypotheses.  By Lemma 2, $2^{(3^{k-1})} \equiv -1 \mod{3^k}$.  

\noindent
$2^{\left(3^{k-1}\right)} \equiv -1 \mod{3^k}$ \\
$ \Rightarrow 2^{3^{\left(k-1\right)}m} \equiv \left( -1\right)^{m} \mod{3^k} $ \\
$\Rightarrow  2^l +1 \equiv \left( -1\right)^{m}+1 \equiv 0  \mod{3^k}$, since $m$ is an odd integer.

\noindent
Because $2^l+1$ is a multiple of $3^k$, it follows that $N$ is an odd integer.

\noindent
\textcolor{blue}{\textbf{Example:}}
$k = 2,\  l = 3^1\cdot 5=15 \to N=\left(\frac{2}{3}\right)^2\left(2^{15}+1\right)-1 = 14563.$
The Collatz sequence for 14563 is:  $ 14563\to21845\to32768\to16384\to 8192\to4096\to \cdots \to8\to4\to2\to1$ (17 steps)

Corollaries 3 and 4 provide special cases of starting values for Collatz sequences that satisfy the Collatz conjecture and for which we can provide the number of steps required to reach 1.
 
\label{Corollary 3}
\noindent
\textcolor {blue}{\textbf{Corollary 3:} } If $q$ is an integer and $N =  \dfrac{4^q-1}{3}$, then the Collatz sequence starting  with $N$ reaches 1 in $2q$ steps.

\noindent
{\textbf{Proof:}  Let $k=1$ in Theorem 2 and simplify algebraically.

\noindent
\textcolor {blue}{\textbf{Example:} } $q=4 \to N=\dfrac{4^4-1}{3} = 85$. \\
\noindent
The Collatz sequence for 85 is: $85 \to128\to 64\to32\to 16\to8\to 4\to2\to 1$. (8 steps).

\label{Corollary 4}
\noindent
\textcolor {blue}{\textbf{Corollary 4:} } If $r$ is an odd integer and $N =  \dfrac{2^{3r+2}-5}{9}$, then the Collatz sequence starting  with $N$ reaches 1 in $3r+2$ steps.

\noindent
\textbf{Proof:}  Let $k=2$ in Theorem 2 and simplify algebraically. 

\noindent
\textcolor {blue}{\textbf{Example:} } $r=3 \to N=\dfrac{2^{11}-5}{9} =227$. \\
\noindent
The Collatz sequence for 227 is:  $227\to341\to512\to 256\to128\to 64\to 32\to16\to 8\to4\to2\to1$.  (11 steps)

We now derive some results which identifies some properties of general Collatz sequences.  We begin by categorizing positive odd integers.

\label{Definition Integer Types}
\noindent 
\textcolor {blue}{\textbf {Definitions: }} Every odd positive integer $N$ falls into one of three categories:
 
\noindent (i)\ \,	$N = 6r + 1,  r = 0,  1, 2,. . ., $ (call this \textbf{type A}) 
Examples:  1, 7, 13,...\\
\noindent (ii)\ \ 	$N = 6r + 3,  r = 0,  1, 2,. . ., $ (call this \textbf{type B})  
Examples:  3, 9, 15,...\\
\noindent (iii)	$N = 6r + 5,  r = 0,  1,  2,.  . ., $ (call this \textbf{type C})  
Examples:  5, 11,  17,...

\label{Theorem 4}
\noindent 
\textcolor {blue}{\textbf {Theorem 4: }} Every odd integer in a Collatz sequence except 
possibly the first one is of type A or C.

\noindent 
\textbf{Proof:}  Let $N$ be an odd integer of type B in a Collatz sequence.  Then there is an integer $r$ such that $N= 6r+3$.  Suppose the number in the Collatz sequence preceding $N$ in the sequence is the odd number 
\noindent $N'$.  Then by the Collatz algorithm, $\dfrac{3N'+1}{2} =N = 6r+3$.  We then have $3N'+1=12r+6$.  Since 3 divides the right-hand side of this equation it must divide the left-hand side, which is impossible.  We conclude that $N$ cannot be immediately preceded by an odd integer.  

Therefore we can assume that either $N$ is at the start of the Collatz sequence or there is an integer $k$ such that the even numbers $2^kN,\  2^{k-1}N, \ \cdots,\  2N$ comprise the sequence before $N$.  This implies that there is an odd number $N'$ preceding $2^kN$.  This means that $\dfrac{3N+1}{2} =2^kN = 2^k(6r+3)=2^k3(2r+1)$. So $3N'+1 = 3(2^{k+1})(2r+1)$.  Since 3 divides the right hand side of this equation, it must divide the left side, which is a contradiction and the theorem is proved.  

\noindent
\textcolor {blue}{\textbf{Remark:}}  This means that no Collatz sequence contains an odd 
multiple of three; i.e., Type B except possibly at the start of the sequence.

\noindent
\textcolor {blue}{\textbf{Example:} } The Collatz sequence for $N = 9$:
$9$(type B)$\to 14\to 7$(type A) $\to 11 $(type C)$\to 17 $(type C)$ \to 26\to 13 $(type A)$ \to 20\to 10\to 5 $(type C) $\to 8\to 4\to 2\to 1$(type A).

The following is a conjecture that expresses a relationship between an odd integer in a Collatz sequence and the next odd integer.

\label{Conjecture}
\noindent
\textcolor {blue}{\textbf{Conjecture:} } Let $N$ be an odd integer in a Collatz sequence and 
let $N' = \tfrac{3N + 1}{2}$ be the next term in the sequence.  If $N'$ factors into $2^lm$ then $m$, the next odd integer in the sequence, is of type A if $l$ is odd and is of type C if $l$ is even.

\noindent
\textcolor {blue}{\textbf{Examples:} }\\
\noindent
$N =15\to N' = 23 = 2^0\cdot23$; $l = 0$ (even) and $m = 23$ (type C).

\noindent
$N =117\to N' = 176 = 2^4\cdot11$; $l = 4$ (even) and $m = 11$ (type C).

\noindent
$N =49\to N' = 74 = 2^1\cdot37$; $l = 1$ (odd) and $m = 37$ (type A).

\noindent
$N =133\to N' = 200 = 2^3\cdot25$; $l = 3$ (odd) and $m = 25$ (type A).

\noindent
$N =341\to N' = 512 = 2^9\cdot1$; $l = 9$ (odd) and $m = 1$ (type A).

\end{document}